\documentclass[leqno]{article}

\def\re{\mathop{\rm Re}}
\def\im{\mathop{\rm Im}}

\def\spec{\mathop{\rm spec}}
\def\supp{\mathop{\rm supp}}

\newtheorem{theorem}{Theorem}
\newtheorem{lemma}[theorem]{Lemma}
\newtheorem{proposition}[theorem]{Proposition}
\newtheorem{definition}[theorem]{Definition}
\newtheorem{corollary}[theorem]{Corollary}

\newcommand{\begintheorem}{\addtocounter{equation}{1}\begin{theorem}}
\newcommand{\beginlemma}{\addtocounter{equation}{1}\begin{lemma}}
\newcommand{\beginproposition}{\addtocounter{equation}{1}\begin{proposition}}
\newcommand{\begindefinition}{\addtocounter{equation}{1}\begin{definition}}
\newcommand{\begincorollary}{\addtocounter{equation}{1}\begin{corollary}}

\begin{document}

\title{Notes on commutative algebra \\
		and harmonic analysis}

\author{Stephen William Semmes	\\
	Rice University		\\
	Houston, Texas}

\date{}

\maketitle

\tableofcontents

\medskip

	In these notes we assume that the reader is familiar with
basic notions from algebra along the lines of groups, rings, fields,
and vector spaces, as well as basic topics in analysis related to
continuity, convergence, compactness, etc.

\section{Commutative groups}
\label{commutative groups}
\setcounter{equation}{0}

	Let $A$ be a commutative group.  Thus $A$ is a nonempty set
with a distinguished element $0$ and which is equipped with a binary
operation $+$ of addition, which is to say that $a + b$ is defined as
an element of $A$ whenever $a, b \in A$, such that addition is
commutative and associative, $0$ is the identity element with respect
to addition, and each element $a$ of $A$ has a unique additive inverse
$-a$ whose sum with $a$ is equal to $0$.  A subset $B$ of $A$ defines
a subgroup of $A$ if $0 \in B$, if the sum of two elements of $B$ is
an element of $B$, and if the additive inverse of an element of $B$ is
an element of $B$.  In other words, $B$ is a commutative group too,
with respect to the operation of addition inherited from $A$.

	A commutative group may also be called an abelian group.  It
may be that a commutative group contains only the zero element, in
which event the group is said to be trivial.  For any commutative
group $A$, $\{0\}$ and $A$ are subgroups of $A$.

	If $A_1$, $A_2$ are commutative groups, then a mapping $\phi$
from $A_1$ to $A_2$ is called a homomorphism if $\phi$ maps the zero
element of $A_1$ to the zero element of $A_2$ and if $\phi$ applied to
a sum of two elements of $A_1$ is equal to the sum in $A_2$ of the
values of $\phi$ at the two specified elements of $A_1$.  The kernel
of a homomorphism $\phi$ from $A_1$ to $A_2$ is the subgroup of $A_1$
consisting of the elements of $A_1$ sent to the zero element of $A_2$
by $\phi$.

	Conversely, suppose that $A$ is a commutative group and that
$B$ is a subgroup of $A$.  One can define a quotient group $A / B$ in
the usual way, basically by identifying $a, a' \in A$ when $a - a' \in
B$.  There is a canonically associated homomorphism from $A$ onto the
quotient group $A / B$ whose kernel is equal to $B$.

	As a basic example we have the commutative group ${\bf Z}$ of
integers with respect to addition.  For each positive integer $m$ we
get a subgroup $m {\bf Z}$ of ${\bf Z}$ consisting of the integer
multiples of $m$.  The quotient ${\bf Z} / m {\bf Z}$, the integers
modulo $m$, is a commutative group with $m$ elements.

	Let $A$ be a commutative group and let $E$ be a nonempty set
which is either finite or countably infinite.  Let us write
$\mathcal{F}(E, A)$ for the space of functions on $E$ with values in
$A$, which is in fact a commutative group with respect to pointwise
addition of functions and with the constant function equal to the zero
element of $A$ at every point in $E$ as the zero element of
$\mathcal{F}(E, A)$.

	If $f$ is a function on $E$ with values in $A$, then the
support of $f$ is denoted $\supp f$ and defined to be the set of $x
\in E$ such that $f(x)$ is not the zero element of $A$.  We define
$\mathcal{F}_0(E, A)$ to be the subgroup of $\mathcal{F}(E, A)$
consisting of functions with finite support.  Of course
$\mathcal{F}_0(E, A)$ is the same as $\mathcal{F}(E, A)$ when $E$ is a
finite set.

\section{Commutative rings}
\label{commutative rings}
\setcounter{equation}{0}

	Let $\mathcal{R}$ be a commutative ring.  Thus $\mathcal{R}$
is a nonempty set with a distinguished element $0$ and two binary
operations of addition and multiplication such that $\mathcal{R}$ is a
commutative group with respect to addition, with $0$ as the additive
identity element, $\mathcal{R}$ is a commutative semigroup with
respect to multiplication, and the two binary operations are
compatible in the sense that usual distributive laws hold.  In
particular, $0$ times any element of the ring is equal to $0$.

	It may be that the ring consists only of $0$ with the rules
that the sum and product of $0$ with itself is equal to $0$.  By a
commutative ring with unity we mean a commutative ring with a nonzero
multiplicative identity element, often denoted $e$ or $1$, which is
then unique.  

	Of course the integers ${\bf Z}$ form a ring using the usual
operations of addition and multiplication.  For each positive integer
$m$, we can consider the integers modulo $m$, which defines a
commutative ring since equivalence modulo $m$ is preserved by both
addition and multiplication.  When $m \ge 2$ this is a commutative
ring with unity.

	Let $E$ be a nonempty set which is either finite or countably
infinite, and let $\mathcal{R}$ be a commutative ring.  As in the
previous section, we can define $\mathcal{F}_0(E, \mathcal{R})$ and
$\mathcal{F}(E, \mathcal{R})$ to be the spaces of functions on $E$
with values in $\mathcal{R}$ with finite support or no restriction on
the support.  Namely, because $\mathcal{R}$ is a commutative group
with respect to addition, these spaces of functions are commutative
groups with respect to pointwise addition as before.  Moreover, we can
use the operation of multiplication on $\mathcal{R}$ to define
multiplication of functions with values in $\mathcal{R}$ pointwise, so
that $\mathcal{F}_0(E, \mathcal{R})$, $\mathcal{F}(E, \mathcal{R})$
become rings.  If $\mathcal{R}$ is a commutative ring with unity, then
so is $\mathcal{F}(E, \mathcal{R})$, with the constant function equal
to the multiplicative identity element of $\mathcal{R}$ as the
multiplicative identity element of $\mathcal{F}(E, \mathcal{R})$.

\section{Fields}
\label{fields}
\setcounter{equation}{0}

	By a field we mean a commutative ring with unity in which
every nonzero element has a multiplicative inverse.  The rational
numbers ${\bf Q}$, the real numbers ${\bf R}$, and the complex numbers
${\bf C}$ are fields.  If $p$ is a prime number, then the integers
modulo $p$ define a field.

	Let $k$ be a field, with multiplicative identity element $1$.
If for each positive integer $n$ the sum of $n$ $1$'s in $k$ is not
equal to $0$ in $k$, then we say that $k$ has characteristic $0$.
Otherwise we say that $k$ has positive characteristic.  If $k$ is a
field with positive characteristic and $n$ is the smallest positive
integer such that the sum of $n$ $1$'s in $k$ is equal to $0$ in $k$,
then $n$ is a prime number, called the characteristic of $k$.  For
each prime number $p$, the integers modulo $p$ form a field with
characteristic $p$.

	Let $k$ be a field.  As usual a vector space over $k$ is a
commutative group equipped with an additional operation of scalar
multiplication by elements of $k$ and with natural compatibility
properties for the field operations on $k$ and the operations of
vector addition and scalar multiplication on the vector space.  If $E$
is a nonempty set which is at most countable, then $\mathcal{F}_0(E,
k)$, $\mathcal{F}(E, k)$ can be viewed as vector spaces over $k$, with
respect to pointwise addition and multiplication by scalars.  If $E$
is finite, so that $\mathcal{F}_0(E, k)$ and $\mathcal{F}(E, k)$ are
the same, then this vector space has dimension equal to the number of
elements of $E$.  Indeed, one can get a basis for the vector space
using functions that are equal to $1$ at one point in $E$ and to $0$
at all other elements of $E$.

	By a commutative algebra over $k$ we mean a vector space over
$k$ equipped with a binary operation of multiplication which is
commutative and associative and satisfies the usual distributive laws
with respect to the vector space operations, which amounts to saying
that multiplication on the algebra is linear in each variable.  One
might also think of a commutative algebra over $k$ as a commutative
ring which happens to be a vector space over $k$ and where
multiplication in the algebra basically commutes with scalar
multiplication in the vector space structure.  If there is a nonzero
multiplicative identity element, then we have a commutative algebra
over $k$ with unity.  For any nonempty set $E$ which is at most
countable the vector spaces $\mathcal{F}_0(E, k)$, $\mathcal{F}(E, k)$
are commutative algebras over $k$ using pointwise multiplication of
functions as multiplication in the algebra.  Also, $\mathcal{F}(E, k)$
is a commutative algebra over $k$ with unity, using the constant
function on $E$ equal to $1$ in $k$ at each point as the
multiplicative identity element.

\section{Convolution algebras}
\label{convolution algebras}
\setcounter{equation}{0}

	Let $E$ be a nonempty set which is finite or countable
equipped with a commutative and associative binary operation $\circ$
and for which there is an element of $E$ which is an identity element
for this operation.  In other words, $E$ is a commutative semigroup
with identity.  Of course a commutative group, like ${\bf Z}$, is a
commutative semigroup with identity element.  The nonnegative integers
or whole numbers ${\bf W}$ form a commutative semigroup under
addition.

	Let $k$ be a field, and consider functions on $E$ with values
in $k$.  If $f_1$, $f_2$ are two such functions, then we would like to
define the convolution $f_1 * f_2$ to be the function on $E$ with
values in $k$ given by
\begin{equation}
	(f_1 * f_2)(z) 
		= \sum_{x, y \in E \atop z = x \circ y} f_1(x) \, f_2(y),
\end{equation}
i.e., where for a given $z \in E$ we sum over all $x, y \in E$ whose
product $x \circ y$ in $E$ is equal to $z$.  If $E$ is finite, then
this certainly makes sense, but if $E$ is infinite, then we should be
more careful.

	If at least one of $f_1$, $f_2$ has finite support, then we
can define the convolution $f_1 * f_2$ as above, since for each $z \in
E$ all but at most finitely many terms in the sum are then equal to
$0$, so that the sum is really a finite sum.  Moreover, if both $f_1$,
$f_2$ have finite support, then the convolution $f_1 * f_2$ has finite
support too.

	Let us define $\mathcal{C}_0(E, k)$ to be the space of
functions on $E$ with values in $k$ and finite support, as a vector
space over $k$ and as an algebra with convolution as multiplication.
In other words, $\mathcal{C}_0(E, k)$ is the same as $\mathcal{F}_0(E,
k)$ as a vector space over $k$, but normally we would use pointwise
multiplication on $\mathcal{F}_0(E, k)$ while now we use convolution.
It is easy to see that convolution satisfies the required properties
in order for $\mathcal{C}_0(E, k)$ to be an algebra over $k$, namely,
it is commutative, associative, and linear in each of the functions
being convolved.

	We can also describe convolution on $\mathcal{C}_0(E, k)$ as
follows.  For each $x \in E$, let $\delta_x(w)$ be the function on $E$
with values in $k$ such that $\delta_x(w)$ is equal to $1 \in k$ when
$w = x$ and to $0 \in k$ when $w \ne x$.  By definition, every element
of $\mathcal{C}_0(E, k)$ is a finite linear combination of the
$\delta_x$'s.  If $x, y \in E$, then one can check that
\begin{equation}
	\delta_x * \delta_y = \delta_{x \circ y}.
\end{equation}
In particular, if $e$ is the identity element of $E$, then $\delta_e$
defines a multiplicative identity element in $\mathcal{C}_0(E, k)$.

	As mentioned earlier, the convolution of two functions on $E$
with values in $k$ makes sense if at least one of the two functions
has finite support.  The convolution of $\delta_e$ with any function
$f$ on $E$ with values in $k$ is equal to $f$.

	In some situations the semigroup $E$ has the nice property
that for each $z \in E$ there are only finitely many pairs of elements
$x$, $y$ in $E$ such that $z = x \circ y$.  For instance, this is the
case for the semigroup of whole numbers.  When this occurs, the
convolution of two functions on $E$ with values in $k$ is always
well-defined, and we define $\mathcal{C}(E, k)$ to be the space of
functions on $E$ with values in $k$ as a vector space and using
convolution as multiplication.  This defines a commutative algebra
over $k$, as one can again easily verify.  Indeed, it is a commutative
algebra with unity because $\delta_e$ defines a nonzero multiplicative
identity element.

\section{Polynomials}
\label{polynomials}
\setcounter{equation}{0}

	Let $k$ be a field, and let $\mathcal{P}(k)$ denote the
algebra of polynomials over $k$.  Thus each element $p(t)$ of
$\mathcal{P}(k)$ can be expressed as
\begin{equation}
	p(t) = c_m \, t^m + c_{m-1} \, t^{m-1} + \cdots + c_0
\end{equation}
for some $c_0, \cdots, c_m \in k$.  One adds and multiplies
polynomials in the usual manner, following the rule that $t^j$ times
$t^l$ is equal to $t^{j + l}$, in such a way that $\mathcal{P}(k)$ is
indeed a commutative algebra over $k$.  It is a commutative algebra
with unity, with the polynomial with only a constant term equal to $1
\in k$ as the multiplicative identity element.

	Let us emphasize that a polynomial here is considered as a
formal expression, with the constant polynomial $1$ and the positive
powers of the indeterminant $t$ being linearly independent by
definition.  In fact one can think of $\mathcal{P}(k)$ as being
another way to write the convolution algebra $\mathcal{C}_0({\bf W},
k)$, where ${\bf W}$ denotes the whole numbers viewed as a semigroup
with respect to addition.  This is equivalent in terms of the algebra
structure, and sometimes nicer as a way to look at the algebra.

	Of course each formal polynomial in $\mathcal{P}(k)$ defines a
polynomial function on $k$.  More precisely, if $p(t) \in
\mathcal{P}(k)$ and $a \in k$, then the value of the function on $k$
associated to $p(t)$ is normally written $p(a)$ and is given by the
same expression as for $p(t)$, with $t$ replaced by $a$.  More
generally, if $\mathcal{A}$ is any algebra over $k$ with unity, then
one gets a function on $\mathcal{A}$ associated to any formal
polynomial $p(t) \in \mathcal{P}(k)$, defined in the same way.  For
example, one can apply this with $\mathcal{A} = \mathcal{P}(k)$, and
this reduces to composing two polynomials following the usual
computation.

	Fix a positive integer $n$, and let us consider polynomials in
$n$ variables.  For this choice of $n$, let us say that a multi-index
is an $n$-tuple $\alpha = (\alpha_1, \ldots, \alpha_n)$ of whole
numbers.  The degree of $\alpha$ is denoted $|\alpha|$ and defined by
\begin{equation}
	|\alpha| = \alpha_1 + \cdots + \alpha_n.
\end{equation}
Using $t = (t_1, \ldots, t_n)$ as an indeterminant, we define the
monomial $t^\alpha$ associated to $\alpha$ by
\begin{equation}
	t^\alpha = t_1^{\alpha_1} \cdots t_n^{\alpha_n}.
\end{equation}
When $\alpha_j = 0$ for some $j$ one can treat $t_j^{\alpha_j}$ as
being equal to $1$.

	By a polynomial in $n$ variables over the field $k$ we mean a
formal linear combination of monomials $t^\alpha$ with coefficients in
$k$.  One can add and multiply polynomials in $n$ variables in the
usual manner, so that the polynomials in $n$ variables over $k$
becomes an algebra over $k$ with unity.  We denote this algebra
$\mathcal{P}_n(k)$, and it can be identified with $\mathcal{C}_0({\bf
W}^n, k)$, where ${\bf W}^n$ denotes the semigroup of $n$-tuples of
whole numbers under addition.  Each formal polynomial in
$\mathcal{P}_n(k)$ determines a polynomial function on $k^n$ in the
usual way, where $k^n$ denotes the space of $n$-tuples with
coordinates in $k$, which can be viewed as a vector space over $k$ of
dimension $n$, and moreover an algebra over $k$ with unity with
respect to coordinatewise multiplication.  More generally, if
$\mathcal{A}$ is any commutative algebra over $k$ with unity, then a
formal polynomial in $\mathcal{P}_n(k)$ leads to a polynomial function
on the space of $n$-tuples with coordinates in $\mathcal{A}$.

	As an instance of this one can take $\mathcal{A}$ to be
$\mathcal{P}_m(k)$ for some positive integer $m$.  This simply
corresponds to the statement that if $p$ is a polynomial of $n$
variables over $k$, and if one has $n$ polynomials $q_1, \ldots, q_n$
of $m$ variables over $k$, then one can get a new polynomial in $m$
variables by replacing the $n$ variables for $p$ with $q_1, \ldots,
q_n$.

	We can be a bit more precise about this by taking homogeneity
into account.  A polynomial $p$ of $n$ variables over $k$ is said to
be homogeneous of degree $l$, where $l$ is a nonnegative integer, if
$p$ can be expressed as a linear combination of monomials $t^\alpha$
with $|\alpha| = l$.  The sum of two homogeneous polynomials of degree
$l$ is homogeneous of degree $l$, and the product of two homogeneous
polynomials of degree $l_1$, $l_2$ is a homogeneous polynomial of
degree $l_1 + l_2$.  A homogeneous polynomial of degree $0$ is the
same as a constant polynomial, and every polynomial can be expressed
in a natural way as a sum of homogeneous polynomials, by grouping
monomials according to their degree.  If $p$ is a homogeneous
polynomial of degree $l$ in $n$ variables over $k$, and if $q_1,
\ldots, q_n$ are homogeneous polynomials of degree $l'$ of $m$
variables over $k$, then when one composes $p$ with $q_1, \ldots, q_n$
as before one gets a homogeneous polynomial of degree $l \cdot l'$ of
$m$ variables over $k$.

	Instead of thinking of polynomials operating on algebras, one
can turn this around a bit and think of elements of algebras as
defining maps on polynomials.  Let $n$ be a positive integer again,
let $\mathcal{A}$ be a commutative algebra with unity over $k$, and
suppose that $a_1, \ldots, a_n$ are elements of $\mathcal{A}$.  Then
the mapping from a polynomial $p$ of $n$ variables over $k$ to its
action on $a_1, \ldots, a_n$ defines a homomorphism from
$\mathcal{P}_n(k)$ into $\mathcal{A}$ as commutative algebras over
$k$.  That is, it is a linear mapping which preserves products, and it
also takes the constant polynomial $1$ to the multiplicative identity
element of $\mathcal{A}$, by definition of the action of a polynomial
on an algebra.  Conversely, every such homomorphism arises in this
manner.

\section{Power series}
\label{power series}
\setcounter{equation}{0}

	Let $k$ be a field, and let $\mathcal{PS}(k)$ denote the
algebra of formal power series over $k$.  In other words, an element
of $\mathcal{PS}(k)$ can be represented as formal sums
\begin{equation}
	\sum_{j=0}^\infty c_j \, t_j,
\end{equation}
where the coefficients $c_j$ are elements of $k$.  One can add and
multiply power series in the usual way, grouping terms corresponding
to the same power of $t$.  As in the case of polynomials, this is
basically equivalent to the convolution algebra $\mathcal{C}(W, k)$
where $W$ is viewed as a semigroup under addition.

	Similarly, we let $\mathcal{PS}_n(k)$ denote the algebra of
formal power series in $n$ variables over $k$.  Elements of
$\mathcal{PS}_n(k)$ can be represented as formal sums
\begin{equation}
	\sum_\alpha c_\alpha \, t^\alpha,
\end{equation}
where the sum extends over all multi-indices $\alpha$, the
coefficients $c_\alpha$ are elements of $k$, and $t^\alpha$ is the
monomial associated to $\alpha$.  Again one can add and multiply power
series in $n$ variables formally, grouping terms which are multiples
of the same monomial, to get an algebra over $k$ with unity.  This
algebra is equivalent to $\mathcal{C}({\bf W}^n, k)$, where ${\bf
W}^n$ is the semigroup of $n$-tuples of whole numbers under addition.

	For each nonnegative integer $l$, let us say that a formal
power series $\sum_\alpha c_\alpha \, t^\alpha$ vanishes to order $l$
if $c_\alpha = 0$ for all multi-indices $\alpha$ such that $|\alpha| <
l$.  In this definition any power series is considered to vanish to
order $0$.  If two power series vanish to order $l$, then so does
their sum, and if two power series vanish to order $l_1$, $l_2$,
respectively, then their product vanishes to order $l_1 + l_2$.  Of
course a power series which vanishes to some order also vanishes to
any smaller order.  A power series that vanishes to all orders is the
zero power series, with all coefficients equal to $0$.

	Let $\sum_\alpha c_\alpha \, t^\alpha$ be a formal power
series in $n$ variables over $k$.  If the constant term in this formal
power series is equal to $0$, then the product of this series with any
other series has constant term equal to $0$ as well.  Conversely, if
the constant term is not equal to $0$, then the series is invertible
in the algebra of formal power series in $n$ variables over $k$, which
is to say that there is another power series whose product with
$\sum_\alpha c_\alpha \, t^\alpha$ is equal to $1$, i.e., the power
series with only a constant term equal to $1$.  Indeed, if $p(t)$ is a
formal power series with constant term equal to $0$, then one can find
an inverse for $1 - p(t)$ as a power series through the formal
expansion $\sum_{l=0}^\infty (p(t))^l$.  A formal power series with
nonzero constant term can be expressed as a nonzero constant times $1
- p(t)$, where $p(t)$ is a power series with constant term equal to
$0$, and hence is invertible.

\section{Rational functions}
\label{Rational functions}
\setcounter{equation}{0}

	Let $k$ be a field, and let us write $\mathcal{R}(k)$ for the
field of rational functions over $k$.  More precisely, a rational
function over $k$ can be represented by a quotient $p(t) / q(t)$,
where $p(t)$, $q(t)$ are polynomials over $k$ and $q(t)$ is not the
zero polynomial, which is to say that it has at least one nonzero
coefficient.  Two such expressions $p_1(t) / q_1(t)$ and $p_2(t) /
q_2(t)$ are viewed as defining the same rational function if $p_1(t)
\, q_2(t) - p_2(t) \, q_1(t)$ is the zero polynomial over $k$.  One
can add and multiply rational functions in the usual manner, and the
reciprocal of a nonzero rational function is a rational function, so
that $\mathcal{R}(k)$ becomes a field which contains a copy of $k$ as
constant functions.

	By a formal Laurent series over $k$ we mean a formal sum of
the form $\sum_{j=a}^\infty c_j \, t^j$, where $a$ is an integer and
the coefficients $c_j$ are elements of $k$.  One can add and multiply
Laurent series as usual by grouping terms which are multiples of the
same power of $t$.  In fact one can take the reciprocal of a nonzero
Laurent series and get a Laurent series through standard computations,
so that the Laurent series form a field which contains $k$ as the
constants.  The rational functions over $k$ form a subfield of the
Laurent series over $k$.

	One can also look at the formal Laurent series in terms of
convolution.  Namely, consider the $k$-valued functions $f(x)$ on the
integers for which there is an integer $a$ such that $f(x) = 0$ when
$x < a$.  These functions form a vector space over $k$ with respect to
pointwise addition and scalar multiplication.  The convolution of two
functions of this type can be defined through the usual formula and
yields another function on ${\bf Z}$ of the same type.  These
functions on ${\bf Z}$ correspond to Laurent series with coefficients
given by the function, with addition and convolution of functions
exactly matching addition and multiplication of Laurent series.

	For each integer $a$, consider the Laurent series which can be
expressed as $\sum_{j = a}^\infty c_j \, t_j$.  Every Laurent series
has this property for some $a$, and if two Laurent series have this
property for a particular $a$ then the sum has the same property.  If
two Laurent series have this property for some $a_1, a_2 \in {\bf Z}$,
then the sum has the same property with $a = a_1 + a_2$.  A Laurent
series that has this property for all integers $a$ is the zero Laurent
series, with all coefficients equal to $0$.

\section{Homomorphisms}
\label{homomorphisms}
\setcounter{equation}{0}

	Let $\mathcal{R}_1$, $\mathcal{R}_2$ be commutative rings.  A
mapping $\phi$ from $\mathcal{R}_1$ into $\mathcal{R}_2$ is said to be
a homomorphism if $\phi$ sends a sum or product of two elements of
$\mathcal{R}_1$ to the corresponding sum or product in
$\mathcal{R}_2$.  If $\mathcal{A}_1$, $\mathcal{A}_2$ are commutative
algebras over the same field $k$, then a homomorphism from
$\mathcal{A}_1$ into $\mathcal{A}_2$ is a ring homomorphism which is
also linear over $k$, so that it cooperates with scalar
multiplication.

	If $\mathcal{R}_1$, $\mathcal{R}_2$ are rings with unity, then
a homomorphism from $\mathcal{R}_1$ into $\mathcal{R}_2$ either maps
all of $\mathcal{R}_1$ to the zero element of $\mathcal{R}_2$, or it
maps the multiplicative identity element of $\mathcal{R}_1$ to the
multiplicative identity element of $\mathcal{R}_2$.  For rings or
algebras with unity one often considers only nonzero homomorphisms.

	Let $\mathcal{R}$ be a commutative ring.  By an ideal in
$\mathcal{R}$ we mean a subset $\mathcal{I}$ of $\mathcal{R}$ which is
a subgroup of $\mathcal{R}$ with respect to addition and which has the
additional property that for each $x \in \mathcal{I}$ and $y \in
\mathcal{R}$ we have that the product $x \, y$ is an element of
$\mathcal{I}$.  If $\mathcal{A}$ is a commutative algebra over a field
$k$, then an ideal in $\mathcal{A}$ should be a linear subspace of
$\mathcal{A}$ as a vector space over $k$ as well as an ideal of
$\mathcal{A}$ as a ring.

	Given a homomorphism between two commutative rings or two
commutative algebras over the same field $k$, the kernel of the
homomorphism is the set of elements of the domain which are sent to
the zero element of the range.  It is easy to check that the kernel of
such a homomorphism is always an ideal.  Conversely, if $\mathcal{I}$
is an ideal in a commutative ring $\mathcal{R}$, then one can form the
quotient $\mathcal{R} / \mathcal{I}$ initially as an abelian group
with respect to addition, and then as a commutative ring, because
multiplication is well-defined on the quotient.  Similarly, if
$\mathcal{A}$ is a commutative algebra over a field $k$, and if
$\mathcal{I}$ is an ideal in $\mathcal{A}$, then one can define the
quotient $\mathcal{A} / \mathcal{I}$ as a vector space over $k$, and
then as an algebra over $k$.  In both cases there is a canonical
quotient homomorphism from the ring or algebra onto the quotient whose
kernel is the ideal $\mathcal{I}$.

	Let us restrict our attention from now on in this section to
commutative rings and algebras with unity.  An ideal in a commutative
ring or algebra is said to be proper if it is not the whole ring or
algebra, which is equivalent to saying that the quotient is not simply
zero in the sense of containing only a zero element.  In a commutative
ring or algebra with unity, an ideal is proper if and only if it does
not contain the multiplicative identity element.  The quotient of a
commutative ring or algebra with unity by a proper ideal is a
commutative ring or algebra with unity, respectively.

	In any commutative ring or algebra one has the zero ideal,
consisting of only the zero element.  This is the only proper ideal in
a field.  If $x$ is an element of a commutative ring or algebra with
unity, then one can define the ideal $(x)$ generated by $x$ to be the
set of all products of the form $x \, y$, where $y$ runs through all
elements of the ring or algebra.  Observe that $(x)$ is equal to the
whole ring or algebra if and only if $x$ is invertible, which is to
say that there is an element $z$ in the ring or algebra such that the
product $x \, z$ is equal to the multiplicative identity element.

	A proper ideal $\mathcal{I}$ in a commutative ring or algebra
with unity is said to be a maximal ideal if the only ideals in the
ring or algebra which contain $\mathcal{I}$ are $\mathcal{I}$ and the
whole ring or algebra.  In the quotient of the ring or algebra by a
maximal ideal $\mathcal{I}$, the only ideals are the zero ideal and
the whole quotient.  This implies that the quotient is a field.

	An ideal $\mathcal{I}$ in a commutative ring or algebra with
unity is said to be a prime ideal if for each $x$, $y$ in the ring or
algebra such that the product $x \, y$ is an element of $\mathcal{I}$
we have that either $x$ or $y$ is an element of $\mathcal{I}$.  In
particular, the zero ideal is prime if and only if there are no
nontrivial zero divisors in the ring or algebra, which is the same as
saying that a product of nonzero elements is nonzero.  The quotient of
a commutative ring or algebra with unity by a prime ideal has no
nontrivial zero divisors, and conversely the kernel of a homomorphism
of a commutative ring or algebra into a commutative ring or algebra,
respectively, with no nontrivial zero divisors is a prime ideal in the
domain.  If a commutative ring or algebra with unity has no nontrivial
zero divisors, then one can form the field of quotients which contains
the original ring or algebra.

\section{Banach spaces}
\label{banach spaces}
\setcounter{equation}{0}

	In this and the next two sections we shall use complex numbers
as scalars.  Recall that a complex number $z$ can be expressed as $x +
y \, i$, where $x$, $y$ are real numbers and $i^2 = -1$.  We call $x$,
$y$ the real and imaginary parts of $z$ and denote them $\re z$, $\im
z$, respectively.

	If $z = x + y \, i$ is a complex number, with $x$, $y$ real
numbers, then the complex conjugate of $z$ is denoted $\overline{z}$
and defined to be $x - y \, i$.  One can check that the complex
conjugate of a sum or product of complex numbers is equal to the
corresponding sum or product of complex conjugates.

	If $x$ is a real number, then the absolute value of $x$ is
denoted $|x|$ and is equal to $x$ when $x \ge 0$ and to $-x$ when $x
\le 0$.  Thus we have $|x + y| \le |x| + |y|$ and $|x \, y| = |x| \,
|y|$ for all real numbers $x$, $y$.  If $z = x + y \, i$ is a complex
number, $x, y \in {\bf R}$, then the modulus of $z$ is denoted $|z|$
and is the nonnegative real number such that $|z|^2 = x^2 + y^2$.
Thus $|z|^2 = z \, \overline{z}$ and if $z$ is a real number then the
modulus of $z$ is the same as the absolute value of $z$.  One can
verify that $|z + w| \le |z| + |w|$ and $|z \, w| = |z| \, |w|$ for
all complex numbers $z$, $w$.

	Let $V$ be a vector space over the complex numbers.  By a norm
on $V$ we mean a function $\|v\|$ defined for all $v \in V$ such that
$\|v\|$ is a nonnegative real number for all $v \in V$ which is
equal to $0$ if and only if $v = 0$,
\begin{equation}
	\|\alpha \, v\| = |\alpha| \, \|v\|
\end{equation}
for all complex numbers $\alpha$ and all $v \in V$, and
\begin{equation}
	\|v + w\| \le \|v\| + \|w\|
\end{equation}
for all $v, w \in V$.  The combination of a complex vector space $V$
and a norm on $V$ is called a complex normed vector space.

	Suppose that $V$ is a complex vector space equipped with a
norm $\|\cdot \|$.  We can use the norm to define a metric on $V$,
namely,
\begin{equation}
	\|v - w\|.
\end{equation}
It is easy to see that this does indeed satisfy the requirements of a
metric, as a consequence of the properties of a norm.  Of course the
complex numbers ${\bf C}$ form a $1$-dimensional complex vector space,
and the modulus $|\cdot |$ of a complex number defines a norm on ${\bf
C}$.

	Let $V$ be a complex vector space with norm $\|\cdot \|$.  If
$\{v_j\}_{j=1}^\infty$ is a sequence of vectors in $V$ and $v$ is a
vector in $V$, then $\{v_j\}_{j=1}^\infty$ converges to $v$ if for
each $\epsilon > 0$ there is a positive integer $L$ such that $\|v_j -
v\| < \epsilon$ for all $j \ge L$.  The limit $v$ of the sequence
$\{v_j\}_{j=1}^\infty$ is unique when it exists, and in this event we
write $v = \lim_{j \to \infty} v_j$.  If $\{v_j\}_{j=1}^\infty$,
$\{w_j\}_{j=1}^\infty$ are sequences of vectors in $V$ which converge
to vectors $v$, $w$, respectively, then $\{v_j + w_j\}_{j=1}^\infty$
converges to $v + w$ in $V$.  Similarly, if
$\{\alpha_j\}_{j=1}^\infty$ is a sequence of complex numbers which
converges to a complex number $\alpha$ and $\{v_j\}_{j=1}^\infty$ is a
sequence of vectors in $V$ which converges to the vector $v$, then
$\{\alpha_j \, v_j\}_{j=1}^\infty$ converges to $\alpha \, v$.

	A sequence of vectors $\{v_j\}_{j=1}^\infty$ in $V$ is said to
be a Cauchy sequence if for each $\epsilon > 0$ there is a positive
integer $L$ such that $\|v_j - v_l\| < \epsilon$ for all $j, l \ge L$.
If a sequence converges, then it is a Cauchy sequence.  If every
Cauchy sequence in $V$ converges, which is the same as saying that the
normed vector space $V$ is complete as a metric space, then we say
that $V$ with the norm $\|\cdot \|$ is a Banach space.

	Suppose that $\{v_j\}_{j=0}^\infty$ is a sequence of vectors
in a normed vector space $V$, and consider the associated infinite
series $\sum_{j=0}^\infty v_j$.  We say that this series converges if
the sequence of partial sums $\sum_{j=0}^n v_j$ converges in $V$.  As
usual a necessary condition for this to happen is that $\lim_{j \to
\infty} v_j = 0$.

	We say that $\sum_{j=0}^\infty v_j$ converges absolutely if
$\sum_{j=0}^\infty \|v_j\|$ converges as an infinite series of
nonnegative real numbers, which is equivalent to saying that the
partial sums $\sum_{j=0}^n \|v_j\|$ are bounded.  If
$\sum_{j=0}^\infty v_j$ converges absolutely, then the partial sums of
$\sum_{j=0}^\infty v_j$ form a Cauchy sequence, by a standard
argument.  In particular it follows that $\sum_{j=0}^\infty v_j$
converges in $V$ if $V$ is a Banach space.  Conversely, if $V$ is a
normed vector space such that every absolutely convergent series in
$V$ converges, then $V$ is a Banach space.  Basically, if
$\{w_l\}_{l=1}^\infty$ is a Cauchy sequence in $V$, then one can first
pass to a subsequence to reduce to the case where $\|w_{l+1} - w_l\|
\le 2^{-l}$ for all $l$, and then express the $w_l$'s as partial sums
of an absolutely convergent series in $V$.

	Let $V$ be a complex vector space, and let $W$ be a linear
subspace of $V$.  We can define the quotient $V / W$ as a complex
vector space in the usual way, identifying elements of $V$ when their
difference lies in $W$.  There is a canonical quotient mapping $q$
which is a linear mapping from $V$ onto $V / W$ whose kernel is equal
to $W$.

	Assume further that $V$ is equipped with a norm $\|\cdot
\|_V$, and that $W$ is closed as a subset of $V$ with respect to the
metric associated to the norm.  We can define a norm $\|\cdot \|_{V /
W}$ on the quotient space by
\begin{equation}
	\|q(v)\|_{V / W} = \inf \{\|v + w\| : w \in W\}
\end{equation}
for all $v \in V$.  It is not too difficult to show that this does
indeed define a norm on $V / W$.  In this regard let us point out that
$\|q(v)\|_{V / W} > 0$ when $q(v) \ne 0$ in $V / W$, which is
equivalent to $v \not\in W$, because $W$ is closed.  For each $v \in
V$ we have that
\begin{equation}
	\|q(v)\|_{V / W} \le \|v\|_V
\end{equation}
by construction.

	If $V$ is a Banach space, then so is the quotient space $V /
W$ equipped with the quotient norm just described.  For if one has a
Cauchy sequence in $V / W$, then one can choose a sequence in $V$
which is mapped to the given sequence in $V / W$ by the quotient
mapping $q$ and which is a Cauchy sequence in $V$.  Alternatively, for
each absolutely convergent series in $V / W$ there is an absolutely
convergent series in $V$ whose terms are sent to the terms of the
series in $V / W$ by the quotient mapping.  In either case, the fact
that $V$ is a Banach space implies that one has convergence there,
which implies convergence in the quotient.

\section{Banach algebras}
\label{banach algebras}
\setcounter{equation}{0}

	Let $\mathcal{A}$ be a commutative algebra with unity over the
complex numbers, where the multiplicative identity element in
$\mathcal{A}$ is denoted $e$.  Suppose that $\|\cdot \|$ is a norm on
$\mathcal{A}$ as a complex vector space.  If
\begin{equation}
	\|x \, y\| \le \|x\| \, \|y\|
\end{equation}
for all $x, y \in \mathcal{A}$ and $\|e\| = 1$, then we say that
$\mathcal{A}$ equipped with this norm is a normed algebra.  If
moreover $\mathcal{A}$ is complete with respect to the metric
associated to the norm, so that $\mathcal{A}$ is a Banach space as a
complex normed vector space, then we say that $\mathcal{A}$ is a
Banach algebra.

	Suppose that $\mathcal{A}$ is a commutative Banach algebra
over the complex numbers with unity and norm $\|\cdot \|$.  If $x \in
\mathcal{A}$, then the spectrum of $x$ is denoted $\spec (x)$ and is
defined to be the set of complex numbers $\lambda$ such that $x -
\lambda \, e$ is \emph{not} invertible.  The set of complex numbers
$\lambda$ such that $x - \lambda \, e$ is invertible is sometimes
called the resolvent set associated to $x$.

	Let $x$ be an element of $\mathcal{A}$ such that $\|x\| < 1$.
The series $\sum_{j=0}^\infty x^j$ is absolutely convergent in
$\mathcal{A}$, since $\|x^j\| \le \|x\|^j$ for all $j$.  Here $x^j$ is
interpreted as being $e$ when $j = 0$, and $\|x\|^j$ is interpreted as
being $1$ when $j = 0$.  As usual,
\begin{equation}
	(e - x) \sum_{j=0}^n x^j = e - x^{n+1}
\end{equation}
for all nonnegative integers $n$.  We conclude that $e - x$ is
invertible in $\mathcal{A}$, with inverse given by $\sum_{j=0}^\infty
x^j$.

	Let $x$ be an element of $\mathcal{A}$ and $\lambda$ be a
complex number such that $|\lambda| > \|x\|$.  We can write $x -
\lambda \, e$ as $-\lambda (e - \lambda^{-1} \, x)$.  It follows that
$x - \lambda \, e$ is invertible, since $\|\lambda^{-1} \, x\| =
|\lambda|^{-1} \, \|x\| < 1$.  In other words, if $\lambda$ is a
complex number which lies in the spectrum of $x$, then
\begin{equation}
	|\lambda| \le \|x\|.
\end{equation}

	Let $x$, $y$ be elements of $\mathcal{A}$ such that $x$ is
invertible.  We can write $y$ as $x + (y - x) = x ( e - x^{-1} (x -
y))$, and it follows that $y$ is invertible if
\begin{equation}
	\|y - x\| < \frac{1}{\|x^{-1}\|}.
\end{equation}
In particular the set of invertible elements of $\mathcal{A}$ is an
open subset of $\mathcal{A}$ with respect to the topology defined by
the metric associated to the norm.  The spectrum of any element of
$\mathcal{A}$ is therefore a closed subset of the complex numbers.

	Let $x$ be an element of $\mathcal{A}$, and assume that the
spectrum of $x$ is empty.  Thus $(x - \lambda \, e)$ is defined for
all $\lambda \in {\bf C}$, and in fact this defines a hoomorphic
function of $\lambda$ with values in $\mathcal{A}$.  Moreover this
holomorphic function tends to $0$ as $|\lambda| \to \infty$.  As in
the scalar case, this implies that the function is constant, and hence
$0$, a contradiction.  In short the spectrum of every element of
$\mathcal{A}$ is nonempty.

	Assume that $\mathcal{A}$ is a field, which is to say that
every nonzero element of $\mathcal{A}$ is invertible.  In this event
the spectrum of each element of $\mathcal{A}$ consists of exactly one
element.  One can show that in fact each element of $\mathcal{A}$ is a
scalar multiple of $e$, so that $\mathcal{A}$ reduces to the complex
numbers.

	Here is a more precise relationship between the norm and
spectrum of an element $x$ of $\mathcal{A}$.  If $\lambda$ is a
complex number such that $x^n - \lambda^n \, e$ is invertible in
$\mathcal{A}$ for some positive integer $n$, then $x - \lambda \, e$
is invertible, because $x^n - \lambda \, e$ can be written as the
product of $x - \lambda \, e$ and a polynomial in $x$.  Hence if
$\lambda$ is a complex number in the spectrum of $x$, then $|\lambda|
\le \|x^n\|^{1/n}$ for all positive integers $n$.  Conversely, if $r$
is a positive real number such that $x - \lambda \, e$ is invertible
in $\mathcal{A}$ for all complex numbers $\lambda$ with $|\lambda| <
r$, then for each positive real number $R > r$ one can show that
$\|x^n\| = O(R^n)$ as $n \to \infty$.  As a result, $\lim_{n \to
\infty} \|x^n\|^{1/n}$ is equal to the maximum of $|\lambda|$ for
$\lambda$ in the spectrum of $x$.

	In the context of Banach algebras it is natural to consider
closed ideals, i.e., ideals which are closed with respect to the
topology determined by the norm.  It is easy to see that the closure
of an ideal is still an ideal, since addition and multiplication are
continuous.  A very nice point is that the closure of a proper ideal
in a commutative Banach algebra with unity is still a proper ideal.
This is because the ball centered at $e$ with radius $1$ is not
contained in a proper ideal, since every element of this ball is
invertible.

	Let $\mathcal{A}$ be a commutative Banach algebra with unity,
and let $\mathcal{I}$ be a proper closed ideal in $\mathcal{A}$.  The
quotient space $\mathcal{A} / \mathcal{I}$ is then a complex
commutative algebra with unity, and it inherits a norm from the one on
$\mathcal{A}$ to become a Banach space.  In fact $\mathcal{A} /
\mathcal{I}$ becomes a commutative Banach algebra with unity in this
way, as one can verify.  Note that the kernel of a continuous
homomorphism from $\mathcal{A}$ into another Banach algebra is a
closed ideal in $\mathcal{A}$.  Conversely, every closed ideal in
$\mathcal{A}$ occurs as the kernel of a continuous homomorphism onto
another Banach algebra, namely, the kernel of the quotient mapping
from $\mathcal{A}$ into $\mathcal{A} / \mathcal{I}$.

	Suppose that $\mathcal{I}$ is a proper ideal in a commutative
Banach algebra $\mathcal{A}$ with unity, and that $\mathcal{I}$ is
maximal, so that the only ideals in $\mathcal{A}$ that contain
$\mathcal{I}$ are $\mathcal{I}$ and $\mathcal{A}$ itself.  This
implies that $\mathcal{I}$ is closed, since the closure of
$\mathcal{I}$ is a closed proper ideal containing $\mathcal{I}$.  The
quotient $\mathcal{A} / \mathcal{I}$ is then a complex commutative
Banach algebra with unity which is a field.  By the result mentioned
earlier $\mathcal{A} / \mathcal{I}$ reduces to the complex numbers.
Thus $\mathcal{I}$ is actually the kernel of a continuous homomorphism
from $\mathcal{A}$ onto the complex numbers.

\section{$\ell^1(E)$}
\label{l^1(E)}
\setcounter{equation}{0}

	Let $E$ be a nonempty set which is finite or countable, and
let $h(x)$ be a real-valued function on $E$ such that $h(x) \ge 0$ for
all $x \in E$.  We say that $h(x)$ is summable on $E$ if the sums of
$h(x)$ over finite subsets of $E$ are bounded.  In this case we define
$\sum_{x \in E} h(x)$ to be the supremum of the sums of $h(x)$ over
finite subsets of $E$.  If $E$ is finite then $h(x)$ is automatically
summable on $E$ and $\sum_{x \in E} h(x)$ is the usual finite sum.

	A complex-valued function $f(x)$ on $E$ is said to be
absolutely summable if $|f(x)|$ is a summable function on $E$.  Again
this is automatic when $E$ is finite.  One can check that the
absolutely summable complex-valued functions on $E$ form a vector
space over the complex numbers, which we denote $\ell^1(E)$.  For such
a function $f(x)$ we put
\begin{equation}
	\|f\|_1 = \|f\|_{\ell^1(E)} = \sum_{x \in E} |f(x)|.
\end{equation}
One can also verify that this defines a norm on $\ell^1(E)$.

	If $f(x)$ is an absolutely summable function on $E$, then one
can define the sum $\sum_{x \in E} f(x)$ as a complex number.  This is
simply a finite sum when $E$ is a finite set.  When $E$ is infinite
there are various ways to look at this.  For instance, one can
enumerate $E$ by the positive integers and treat $\sum_{x \in E} f(x)$
as an ordinary infinite series, which is absolutely convergent.  The
value of the sum does not depend on the way that $E$ is enumerated by
the positive integers, because of absolute convergence.

	In fact $\ell^1(E)$ is a Banach space, which is to say that if
$\{f_j\}_{j=1}^\infty$ is a Cauchy sequence in $\ell^1(E)$ with
respect to the norm $\|\cdot \|_1$, then there is an $f \in \ell^1(E)$
such that $\{f_j\}_{j=1}^\infty$ converges to $f$.  Indeed, if
$\{f_j\}_{j=1}^\infty$ is a Cauchy sequence in $\ell^1(E)$, then for
each $x \in E$ we have that $\{f_j(x)\}_{j=1}^\infty$ is a Cauchy
sequence of complex numbers.  It follows that
$\{f_j(x)\}_{j=1}^\infty$ converges as a sequence of complex numbers
for each $x \in E$, and we put $f(x) = \lim_{j \to \infty} f_j(x)$.
One can show that $f \in \ell^1(E)$ by bounding the sums of $|f(x)|$
over finite subsets of $E$ in terms of an upper bound for $\|f_j\|_1$.
The last step is to check that $\{f_j\}_{j=1}^\infty$ converges to $f$
in the $\ell^1(E)$ norm.

	Of course complex-valued functions on $E$ with finite support
are absolutely summable.  The functions on $E$ with finite support
form a linear subspace of $\ell^1(E)$ which is dense with respect to
the topology defined by the norm.  This gives another approach to the
definition of $\sum_{x \in E} f(x)$ for $f \in \ell^1(E)$.  Namely,
$\sum_{x \in E} f(x)$ reduces to a finite sum when $f$ has finite
support, and the general case can be handled by approximation.  For
this we should notice that the modulus of $\sum_{x \in E} f(x)$ is
less than or equal to $\|f\|_1$.

	Now assume that $E$ is a commutative semigroup with identity.
If $f_1$, $f_2$ are two absolutely summable functions on $E$, then we
can define their convolution $f_1 * f_2$ as a function on $E$ as
before.  One can check that $f_1 * f_2$ is also absolutely summable,
and that
\begin{equation}
	\|f_1 * f_2\|_1 \le \|f_1\|_1 \, \|f_2\|_1.
\end{equation}
Recall that the function on $E$ which is equal to $1$ at the identity
element of $E$ and equal to $0$ at all other elements of $E$ is an
identity element for convolution.  It follows that $\ell^1(E)$ is a
commutative Banach algebra with unity over the complex numbers, using
convolution as multiplication and $\|f\|_1$ as the norm.

	Let us consider the special case where $E$ is the group of
integers under addition.  Let ${\bf T}$ denote the set of complex
numbers $z$ such that $|z| = 1$, and let $C(T)$ denote the algebra of
continuous complex-valued functions on ${\bf T}$.  If $\phi(z)$ is an
element of $C({\bf T})$, then $\phi(z)$ is bounded, and we can put
\begin{equation}
	\|\phi \|_{sup} = \sup \{|\phi(z)| : z \in {\bf T}\}.
\end{equation}
This is called the supremum norm of $\phi$, and it is well-known that
$C(T)$ becomes a Banach space with respect to this norm.  Moreover
$C({\bf T})$ is a commutative Banach algebra with unity over the
complex numbers with respect to this norm.

	If $f \in \ell^1({\bf Z})$, then we get an associated function
$\phi_f$ on the unit circle by
\begin{equation}
	\phi_f(z) = \sum_{j=-\infty}^\infty f(j) \, z^j.
\end{equation}
The absolute summability of $f$ ensures that these series of complex
numbers converge absolutely, and it follows from the Weierstrass
$M$-test that the partial sums of the series converge uniformly to a
continuous function on ${\bf T}$.  Clearly
\begin{equation}
	\|\phi_f\|_{sup} \le \|f\|_1
\end{equation}
for all $f \in \ell^1(E)$, and the mapping $f \mapsto \phi_f$ is a
continuous homomorphism from $\ell^1(E)$ as a Banach algebra using
convolution into $C({\bf T})$ as a Banach algebra with respect to
pointwise multiplication of functions.  This mapping is also
one-to-one.

	Recall that ${\bf W}$ denotes the semigroup of whole numbers
under addition, which is a sub-semigroup of the group of integers
under addition.  We can identify $\ell^1(W)$ with a linear subspace of
$\ell^1({\bf Z})$, namely with the subspace of functions $f(j)$ such
that $f(j) = 0$ when $j < 0$, and this subspace is a subalgebra of
$\ell^1(E)$ with respect to convolution.  For $f \in \ell^1({\bf W})$
we can define $\phi_f(z)$ as a continuous function on the closed unit
disk $\{z \in {\bf C} : |z| \le 1\}$ using exactly the same formula as
before.  Again we have that $|\phi_f(z)| \le \|f\|_1$ for all $z \in
{\bf C}$ with $|z| \le 1$ and all $f \in \ell^1(E)$, and we can think
of $f \mapsto \phi_f$ as a continuous homomorphism from $\ell^1(W)$
into the commutative Banach algebra with unity of continuous
complex-valued functions on the unit disk which are complex-analytic
on the open unit disk.

\end{document}